\documentclass[final]{l4dc2024pp}

\title[Piecewise regression via mixed-integer programming for MPC]{Piecewise regression via mixed-integer programming for MPC}
\usepackage{times}
\usepackage{booktabs}
\usepackage[nice]{nicefrac}

\author{%
 \Name{Dieter Teichrib} \Email{dieter.teichrib@tu-dortmund.de}\\
 \Name{Moritz {Schulze~Darup}} \Email{moritz.schulzedarup@tu-dortmund.de}\\
 \addr Control and Cyberphysical Systems Group, TU Dortmund University, Germany%
}

\newcommand{\R}{\mathbb{R}}

\newcommand{\Cc}{\mathcal{C}}
\newcommand{\Dc}{\mathcal{D}}

\newcommand{\Ic}{\mathcal{I}}

\newcommand{\Rc}{\mathcal{R}}

\newcommand{\Wb}{\boldsymbol{W}}
\newcommand{\Vb}{\boldsymbol{V}}

\newcommand{\ab}{\boldsymbol{a}}
\newcommand{\bb}{\boldsymbol{b}}

\newcommand{\gb}{\boldsymbol{g}}
\newcommand{\hb}{\boldsymbol{h}}

\newcommand{\xb}{\boldsymbol{x}}

\newcommand{\oneb}{\boldsymbol{1}}
\newcommand{\zerob}{\boldsymbol{0}}

\newcommand{\gammab}{\boldsymbol{\gamma}}
\newcommand{\deltab}{\boldsymbol{\delta}}

\newcommand{\thetab}{\boldsymbol{\theta}}

\begin{document}

\maketitle

\begin{abstract}%
Piecewise regression is a versatile approach used in various disciplines to approximate complex functions from limited, potentially noisy data points. In control, piecewise regression is, e.g., used to approximate the optimal control law of model predictive control (MPC), the optimal value function, or unknown system dynamics. Neural networks are a common choice to solve the piecewise regression problem. However, due to their nonlinear structure, training is often based on gradient-based methods, which may fail to find a global optimum or even a solution that leads to a small approximation error. To overcome this problem and to find a global optimal solution, methods based on mixed-integer programming (MIP) can be used. However, the known MIP-based methods are either limited to a special class of functions, e.g., convex piecewise affine functions, or they lead to complex approximations in terms of the number of regions of the piecewise defined function. Both complicate a usage in the framework of control. We propose a new MIP-based method that is not restricted to a particular class of piecewise defined functions and leads to functions that are fast to evaluate and can be used within an optimization problem, making them well suited for use in control.  
\end{abstract}

\vspace{-2mm}
\section{Introduction}

Piecewise regression is used in many disciplines to fit a continuous piecewise defined function to an unknown function based on a finite number of possibly noisy samples. In the context of model predictive control (MPC) for linear systems, piecewise regression is used, e.g., to approximate the piecewise affine (PWA) optimal control law \citep{Chen2018,Karg2020} or the piecewise quadratic optimal value function (OVF) \citep{Teichrib2021}. Both functions are continuous and defined on polyhedral regions, which partition the state space (c.f. \citet{Bemporad2002}). In addition, piecewise regression may be used to find a surrogate model for an unknown nonlinear system that is used as a prediction model in nonlinear MPC. Neural networks (NN) are very popular in piecewise regression as they allow to approximate a large class of functions with arbitrary accuracy \citep{Hornik1989}. However, due to the nonlinear structure of NN, computing NN parameters, which minimize the approximation error, requires solving a nonlinear optimization problem (OP). This OP is typically solved approximately by gradient-based methods, such as stochastic gradient descent, RMSProb, or Adam \citep{kingma2017adam}. However, these methods are often sensitive to the choice of hyperparameters, and they may even fail to converge for some choices of the learning rate. Moreover, even if the hyperparameters are chosen properly, general convergence guarantees are lacking for the popular RMSProb and Adam algorithms. Furthermore, in case of approximated system dynamics or OVF, trained NN may subsequently be used as a prediction model or to derive optimal control actions, respectively. For both applications, the NN will be included in an OP (either as constraints or costs), and an efficient inclusion is hard for some types of NN. We will consider maxout NN with a special topology that allows to solve the problem of piecewise regression globally and that is suitable for a subsequent usage in an OP.

An alternative for solving the problem of piecewise regression is the direct use of PWA functions. When using PWA functions for regression, either the regions are fixed and a linear regression is performed for each region, or the computation is split into two parts. In the first part, the data is divided into clusters that define the polyhedral regions of the PWA function. In the second part, the parameters of the function are computed by performing a linear regression for every region (see, e.g., \citet{Magnani2009}). Often, these algorithms alternate between clustering and regression until a termination criterion is met \citep{Bemporad2023}. The main drawback of these methods is that the regression problem is only solved approximately and the approximation quality depends on the initial clustering of the data points. To overcome this problem and to find a global optimal solution, clustering and regression should be done simultaneously. Most of the approaches where clustering and regression are formulated in a combined OP are based on mixed-integer programming (MIP). One of the first MIP-based approaches was introduced in \citet{Roll2004}, where so-called hinging hyperplanes \citep{Breiman1993} are fitted via MIP to the data. However, only PWA functions that possess the consistent variation property \citep[Def.~5]{Chua1988} can be exactly represented by hinging hyperplanes. Since the consistent variation property typically does not hold in MPC, hinging hyperplanes are not the ideal choice. Another MIP-based method was introduced in \citet{Bertsimas2007}. They use binary variables to assign points to regions, we will use similar ideas in our method. In addition, we will also use clustering as in \citet[Sec.~3.2.]{Bertsimas2007} to cluster points and then assign these clusters to regions of the function. The main weakness of this method is that the function is discontinuous and it is not possible to enforce continuity for problems with $n>1$ (c.f. \citet[Sec.~5.3.]{Bertsimas2007}). Moreover, the regions of the function are not necessarily polytopes. Hence, further computations are required to find polytopes that contain the majority of the points associated with an individual region. Furhter MIP-based methods were recently introduced in \citet{Rebennack2020,Warwicker2023}. However, these methods can only be applied to univariate or bivariate functions.
The MIP-based method from \citet[Sec.~3.2.]{Toriello2012} solves some of the aforementioned problems. This method allows the fitting of a convex continuous PWA function to data via a mixed-integer quadratic program (MIQP). The main idea is to represent the convex function as the maximum of its affine segments and use binary variables to determine which affine segment provides the largest value for a given data point. The advantage of this method is that the polyhedral regions of the function are implicitly given by the structure of the function and do not need to be computed separately. However, the method is not extendable to non-convex functions. In \citet{Siahkamari2020}, an alternative method is given, which allows to combine the benefits of \citet{Bertsimas2007} and \citet{Toriello2012}, i.e., it is possible to fit non-convex continuous functions that are defined on polyhedral regions. In \citet{Siahkamari2020}, the fitted PWA function has one region per data point, this allows to solve the piecewise regression problem by solving a quadratic program (QP) instead of an MIQP as in \citet{Roll2004,Bertsimas2007,Toriello2012}. Unfortunately, the coupling between the number of regions and the number of data points leads to rather complex functions for large data sets, which is disadvantageous, especially when the function is used in an OP.

\textbf{Contributions.} The main contribution of this paper is to provide a new MIP-based method to globally solve the piecewise regression problem, which combines the advantages of the known approaches. More specifically, our new method 

1. \quad allows to freely choose the number of segments as in \citet{Bertsimas2007} without the need to compute the regions of the function explicitly,

2. \quad has a simple representation of the fitted function as in \citet{Toriello2012} without being restricted to convex functions,
    
3. \quad allows to fit a continuous function as in \citet{Siahkamari2020} without coupling the number of regions and the number of data points,

4. \quad results in a model that is suitable for the use in the constraints and the cost function of an OP as in \citet{Bemporad2023}.

We also provide two approaches for reducing the complexity of the proposed method to make it applicable to larger data sets. 
        
\vspace{-2mm}
\section{Problem formulation and preliminaries}\label{sec:Preliminaries}
We consider the problem of continuous piecewise regression, i.e., our aim is to fit a continuous piecewise defined function $\Phi(\xb;\thetab):\R^n\rightarrow\R$ to samples of an unknown function $f(\xb):\R^n\rightarrow\R$, which typically involves solving the nonlinear OP 
\vspace{-1mm}
\begin{equation}
    \label{eq:trainigNonLinear}
    \min\limits_{\thetab}\sum_{k=1}^{N_D} \ell \big( y^{(k)}-\Phi(\xb^{(k)};\thetab) \big),
    \vspace{-1mm}
\end{equation}
for a given data set $\Dc:=\{ (\xb^{(i)},y^{(i)}) \in \Rc^{n+1} \ | \ y^{(i)}=f(\xb^{(i)}), i \in \{1,\dots,N_D\} \}$, where $\thetab$ are the parameters of the function and $\ell$ is a convex loss function. In the following, we briefly summarize the methods from \citet{Toriello2012} and \citet{Siahkamari2020}, as our method is most closely related to these and contains both as special cases. In \citet{Toriello2012}, the problem of piecewise affine regression is solved by the MIQP
\vspace{-1mm}
\begin{subequations}\label{eq:OP_trainingCPWA_MIQP}
\begin{align}
    \min_{\Vb,\alpha^{(i)},\deltab^{(i)}} \sum_{i=1}^{N_D} \big( y^{(i)}&-\alpha^{(i)} \big)^2 \\
    \label{eq:OPCPWA_V}
    \text{s.t.} \hspace{10mm} \Vb \begin{pmatrix} {\xb^{(i)}} \\ 1 \end{pmatrix} &\leq \oneb \alpha^{(i)}, \quad \Vb \begin{pmatrix} {\xb^{(i)}} \\ 1 \end{pmatrix} + M (\oneb-\deltab^{(i)}) \geq \oneb \alpha^{(i)}, \quad \forall i \in \{1,\dots,N_D\}, \\
    \sum\limits_{k=1}^{p_1} \deltab_k^{(i)} &= 1, \ \deltab^{(i)} \in \{0,1\}^{p_1}, \hspace{44mm} \forall i \in \{1,\dots,N_D\}. 
    \vspace{-1mm}
\end{align}
\end{subequations}
with $\Vb\in\R^{p_1 \times (n+1)}$ and the subscript $j$ refers to the $j$-th row of a matrix or column vector. The $N_D p_1$ binary variables select which of the $p_1$ affine segments $\Vb_j ( {\xb^{(i)}}^\top \quad 1 )^\top$ with $j\in\{1,\dots,p_1\}$ provides the largest value for a given data point. The constraints~\eqref{eq:OPCPWA_V} ensure that $\max\{\Vb \gb(\xb^{(i)})\}=\alpha^{(i)}$ holds, i.e., for a sufficiently large $M\in\R_{>0}$, the MIQP~\eqref{eq:OP_trainingCPWA_MIQP} solves the nonlinear OP~\eqref{eq:trainigNonLinear} with $\ell(\cdot):=\lVert \cdot \rVert_2^2$, $\thetab:=\Vb$, and $\Phi(\xb;\thetab)=\max\left\{\Vb \begin{pmatrix} \xb \\ 1 \end{pmatrix} \right\}:=\max\limits_{1\leq i \leq p_1}\left\{\Vb_i \begin{pmatrix} \xb \\ 1 \end{pmatrix} \right\}$. Thus, a convex PWA function is fitted to the data set. In contrast to \citet{Bertsimas2007}, this method allows to make predictions by evaluating the function $\Phi(\xb)$ without explicitly computing the regions of the function. Moreover, the function $\Phi(\xb)$ is continuous by construction. However, the disadvantage is that the fitted function is always convex and PWA. In our method, we also use the binary variables to select the largest segments, but we will overcome the limitation to convex functions by extending the MIQP~\eqref{eq:OP_trainingCPWA_MIQP}.

A method that combines the ability to fit non-convex functions with the continuity property of \eqref{eq:OP_trainingCPWA_MIQP} is presented in \citet{Siahkamari2020}, where the problem of piecewise affine regression is solved by solving a QP (see \citet[Eq.~(6)]{Siahkamari2020}) instead of an MIQP. The constraints in the QP ensure that the function
\vspace{-1mm}
\begin{equation}
    \label{eq:PhiSia}
    \Phi(\xb)=\max_{1\leq i \leq N_D}\{ \ab_i^\top (\xb-\xb^{(i)}) + \hat{y}_i + z_i \} - \max_{1\leq i \leq N_D}\{ \bb_i^\top (\xb-\xb^{(i)}) + z_i \}
    \vspace{-1mm}
\end{equation}
is such that we have $\Phi(\xb^{(i)})=\hat{y}_i$ and the cost function is the sum of the residuals $(y^{(i)}-\hat{y}_i)^2$, i.e., the function \eqref{eq:PhiSia} is fitted to the data. Since every PWA can be decomposed into the difference of two convex PWA functions and \eqref{eq:PhiSia} is the difference of two convex PWA functions, a general PWA function is fitted to the data points in \citet{Siahkamari2020}. Unfortunately, unlike \citet{Bertsimas2007} or \citet{Toriello2012}, this method does not allow to choose the number of regions of the PWA function freely. As apparent from \eqref{eq:PhiSia}, the number of regions is always equal to the number of data points, resulting in a rather complex function for large data sets. 

\vspace{-2mm}
\section{A new MIP-based method for piecewise regression}

\begin{figure}[h]
\begin{minipage}[h]{0.51\columnwidth}
\centering
    \subfigure[\hspace{0mm}]{\includegraphics[trim={8.1cm 12.5cm 8.3cm 12.5cm},clip,scale=0.5]{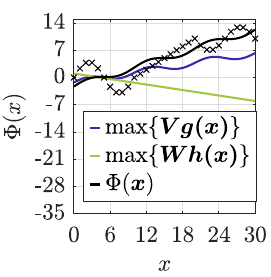}}
    \subfigure[\hspace{0mm}]{\includegraphics[trim={8.1cm 12.5cm 8.3cm 12.5cm},clip,scale=0.5]{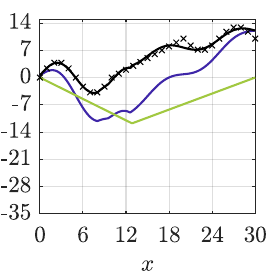}}
    \subfigure[\hspace{0mm}]{\includegraphics[trim={8.1cm 12.5cm 8.3cm 12.5cm},clip,scale=0.5]{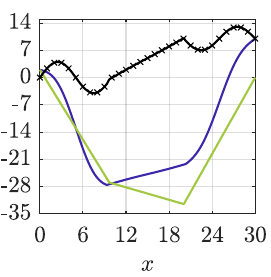}}
\caption{$\Phi(x)$ with $p_i=1$ in (a), $p_i=2$ in (b), and $p_i=3$ in (c) for $i\in\{1,2\}$.}\label{fig:PWS}
\end{minipage}\quad
\begin{minipage}[h]{0.51\columnwidth}
    \subfigure[\hspace{0mm}]{\includegraphics[trim={8.1cm 12.5cm 8.3cm 12.5cm},clip,scale=0.5]{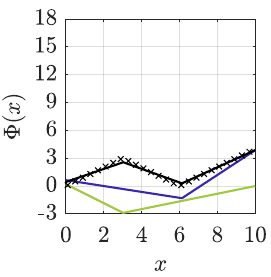}}
    \subfigure[\hspace{0mm}]{\includegraphics[trim={8.1cm 12.5cm 8.3cm 12.5cm},clip,scale=0.5]{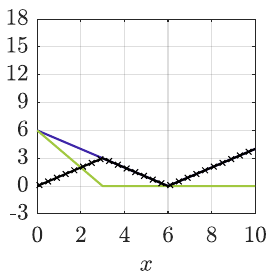}}
    \subfigure[\hspace{0mm}]{\includegraphics[trim={8.1cm 12.5cm 8.3cm 12.5cm},clip,scale=0.5]{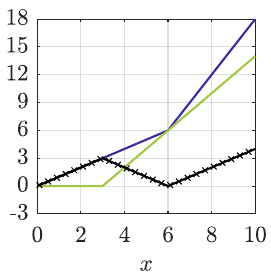}}
\caption{PWA $\Phi(x)$ with $M\!=\!10$ in (a),\\ $M\!=\!50$ in (b), and $M\!=\!100$ in (c).}\label{fig:PWA}
\end{minipage}
\vspace{-2mm}
\end{figure}

In this section, we present our new method for piecewise regression via MIP, which combines advantageous properties of the known approaches while avoiding some of their limitations. Our starting point is the continuous function $\Phi(\xb):\R^n\rightarrow\R$ parameterized by
\vspace{-1mm}
\begin{equation}\label{eq:Phi}
    \Phi(\xb):=\max \{\Vb \gb(\xb) \} - \max \{\Wb \hb(\xb) \}
    \vspace{-1mm}
\end{equation}
with $\Vb\in\R^{p_1 \times r_1}$ and $\Wb\in\R^{p_2 \times r_2}$ as well as $\gb(\xb): \R^n\rightarrow\R^{r_1}$ and $\hb(\xb): \R^n\rightarrow\R^{r_2}$. We denote $\Vb_j \gb(\xb)$ with $j\in\{1,\dots,p_1\}$ and $\Wb_j \hb(\xb)$ with $j\in\{1,\dots,p_2\}$ as segments of the max-functions. In addition, we refer to a segment $i$ as being active if it provides the largest value for a given data point $\xb$, i.e., if $\Vb_i \gb(\xb) \geq \Vb_j \gb(\xb)$ holds for all $j\in\{1,\dots,p_1\}$. This form of parameterization is chosen because every PWA function can be parameterized in this way if $\gb(\xb)$ and $\hb(\xb)$ contain affine and constant functions, i.e., for
\vspace{-1mm}
\begin{equation}\label{eq:ghForPWA}
    \gb(\xb)=\hb(\xb)=
    \begin{pmatrix}
        \xb \\ 1
    \end{pmatrix}.
    \vspace{-1mm}
\end{equation}
Moreover, in \citet{Teichrib2022}, it is shown that every PWQ with $n=1$ can be represented in the form~\eqref{eq:Phi} if $\gb(\xb)$ further contains all quadratic monomials of $\xb$. Alternatively, the function~\eqref{eq:Phi} can be interpreted as a shallow maxout NN \citep{Goodfellow2013} with two neurons of rank $p_1$ and $p_2$, respectively. For the training of such NN, the nonlinear OP~\eqref{eq:trainigNonLinear} with $\Phi(\xb)$ as in \eqref{eq:Phi} is often solved with gradient-based methods. However, these gradient-based methods may get stuck in poor local optima. This problem can be avoided by reformulating the nonlinear OP as the MIQP
\vspace{-1mm}
\begin{subequations}\label{eq:OP_training_MIQP}
\begin{align}
    \min_{\Vb,\Wb,\alpha^{(i)},\beta^{(i)},\gammab^{(i)},\deltab^{(i)}} \sum_{i=1}^{N_D} \big( y^{(i)}-&(\alpha^{(i)} -\beta^{(i)}) \big)^2 \\
    \label{eq:VgLeqAlpha}
    \text{s.t.} \hspace{22.5mm} \Vb \gb(\xb^{(i)}) &\leq \oneb \alpha^{(i)}, \quad \Vb \gb(\xb^{(i)}) + M (\oneb-\deltab^{(i)}) \geq \oneb \alpha^{(i)}, \\
    \label{eq:sumDelta}
    \sum\limits_{k=1}^{p_1} \deltab_k^{(i)} &= 1, \quad \deltab^{(i)} \in \{0,1\}^{p_1}, \\
    \label{eq:WhLeqBeta}
    \Wb \hb(\xb^{(i)}) &\leq \oneb \beta^{(i)}, \quad \Wb \hb(\xb^{(i)}) + M (\oneb-\gammab^{(i)}) \geq \oneb \beta^{(i)}, \\
    \label{eq:sumGamma}
    \sum\limits_{k=1}^{p_2} \gammab_k^{(i)} &= 1, \quad \gammab^{(i)} \in \{0,1\}^{p_2},
    \vspace{-1mm}
\end{align}
\end{subequations}
for all $i \in \{1,\dots N_D\}$. The binary variable vectors $\deltab^{(i)}$ and $\gammab^{(i)}$ in \eqref{eq:VgLeqAlpha} and \eqref{eq:WhLeqBeta} select which segments of the max-functions are active for the point $\xb^{(i)}$. The constant $M\in\R_{>0}$, often referred to as big-$M$, is typically very large. We assume without loss of generality that $p_j\leq N_D$ for $j\in\{1,2\}$, because otherwise there are segments which are inactive for all $\xb^{(i)}$ with $i\in\{1,\dots,N_D\}$. Note that for the special case $p_2=0$ as well as $\gb(\xb)$ and $\hb(\xb)$ as in \eqref{eq:ghForPWA} the MIQP~\eqref{eq:OP_training_MIQP} becomes equivalent to \eqref{eq:OP_trainingCPWA_MIQP}. In the following theorem, we will show that the MIQP~\eqref{eq:OP_training_MIQP} can be used to solve the OP~\eqref{eq:trainigNonLinear}.  
\begin{theorem}\label{thm:TrainingMIQP}
    Let $\Phi(\xb)$ be as in \eqref{eq:Phi}, $\ell(\cdot):=\lVert \cdot \rVert_2^2$, and $\thetab:=\{\Vb,\Wb\}$. Then there exists an $M\in \R_{>0}$ for which the MIQP~\eqref{eq:OP_training_MIQP} is equivalent to the nonlinear OP~\eqref{eq:trainigNonLinear}.
\end{theorem}
\begin{proof}
    The constraints \eqref{eq:sumDelta} are such that for every $i \in \{1,\dots,N_D\}$ there exist exactly one $k_i$ with $\deltab^{(i)}_{k_i}=1$ and $\deltab^{(i)}_{k}=0$ for all $k\in\{1,\dots,p_1\}\setminus k_i$. Thus, the constraints \eqref{eq:VgLeqAlpha} and \eqref{eq:sumDelta} can be equivalently written as
    \vspace{-1mm}
    \begin{subequations}
    \begin{align}
        \label{eq:VEqAlpha}
        &\Vb_{k_i}\gb(\xb^{(i)}) = \alpha^{(i)}, \quad \Vb_k \gb(\xb^{(i)}) \leq \Vb_{k_i} \gb(\xb^{(i)}), \ \forall i \in \{1,\dots,N_D\}, \ \forall k \in \{1,\dots,p_1\} \\
        \label{eq:ViGeqVki}
        &\Vb_k \gb(\xb^{(i)}) + M \geq \Vb_{k_i}\gb(\xb^{(i)}), \hspace{26.5mm} \forall i \in \{1,\dots,N_D\}, \ \forall k \in \{1,\dots,p_1\}.
        \vspace{-1mm}
    \end{align}
    \end{subequations}
    As a consequence, the constraints \eqref{eq:VEqAlpha} are equivalent to $\alpha^{(i)}=\max\{ \Vb \gb(\xb^{(i)}) \}$. By using similar arguments, we can show that similar conditions hold for $\beta^{(i)}$, which leads to the OP
    \vspace{-1mm}
    \begin{subequations}\label{eq:OP_training_MIQP_proof}
    \begin{align}
    \min_{\Vb,\Wb} \sum_{i=1}^{N_D} \big( y^{(i)}-&(\max\{ \Vb \gb(\xb^{(i)}) \} - \max\{ \Wb \hb(\xb^{(i)}) \}) \big)^2 \\
    \label{eq:alphaMax_proof}
    \text{s.t.} \quad \max\{ \Vb \gb(\xb^{(i)}) \} &=\alpha^{(i)}, \hspace{6mm} \max\{ \Wb \hb(\xb^{(i)}) \} =\beta^{(i)}, \hspace{6.1mm} \forall i \in \{1,\dots,N_D\}, \\ 
    \label{eq:VgLeqAlpha_proof}
    \Vb \gb(\xb^{(i)}) + \oneb M &\geq \oneb \alpha^{(i)}, \quad \Wb \hb(\xb^{(i)}) + \oneb M \geq \oneb \beta^{(i)}, \quad \forall i \in \{1,\dots,N_D\}. 
    \vspace{-1mm}
    \end{align}
    \end{subequations}
    This OP is always feasible since the trivial solution ($\Vb=\zerob$, $\Wb=\zerob$, and $\alpha^{(i)}=\beta^{(i)}=0$ for all $i \in \{1,\dots,N_D\}$) is in the feasible set. Furthermore, since $\oneb\alpha^{(i)}-\Vb\gb(\xb^{(i)})$ and $\oneb\beta^{(i)}-\Wb\hb(\xb^{(i)})$ are finite, we can always find a finite $M\in\R_{>0}$ such that the constraints \eqref{eq:VgLeqAlpha_proof} never hold with equality. For such an $M$, the constraints \eqref{eq:alphaMax_proof} and \eqref{eq:VgLeqAlpha_proof} can be neglected. Consequently, the OP~\eqref{eq:OP_training_MIQP_proof} and thus \eqref{eq:OP_training_MIQP} are equivalent to the nonlinear OP \eqref{eq:trainigNonLinear} with $\Phi(\xb)$ as in \eqref{eq:Phi}, $\ell(\cdot):=\lVert \cdot \rVert_2^2$, and $\thetab:=\{\Vb,\Wb\}$.
\end{proof}
According to Theorem~\ref{thm:TrainingMIQP}, the MIQP~\eqref{eq:OP_training_MIQP} allows to compute a global optimal solution for a special instance of the nonlinear OP~\eqref{eq:trainigNonLinear}. Remarkably, the MIQP~\eqref{eq:OP_training_MIQP} can be used to find an approximation of the form \eqref{eq:Phi} for every continuous $\gb_i(\xb)$ and $\hb_j(\xb)$, e.g., by choosing $\gb(\xb)$ quadratic and $\hb(\xb)$ affine, the function represented by \eqref{eq:Phi} is a PWQ function. Even more complex functions are possible. For example, Figure~\ref{fig:PWS} shows a piecewise defined function with sinusoidal segments, approximated by \eqref{eq:Phi} with $\gb(x)=\begin{pmatrix} 1 & x & \sin(0.2x) \end{pmatrix}^\top$ and $\hb(x)=\begin{pmatrix} 1 & x \end{pmatrix}^\top$. In all figures the functions $\max\{\Vb \gb(\xb) \}$ and $\max\{\Wb \hb(\xb) \}$ are represented in blue and green, respectively. The approximation of $f(x)$ by $\Phi(x)$ is illustrated in black, where the data points $(x^{(i)},f(x^{(i)}))$ used in the MIQP~\eqref{eq:OP_training_MIQP} are the black crosses.

Regarding the $M$ from the MIQP~\eqref{eq:OP_training_MIQP}, it should be noted that methods from the literature (see, e.g., \citet{Bonami2015}) to compute the smallest possible value of $M$ cannot be applied here, because the weights $\Vb$ and $\Wb$ are optimization variables and thus not known in advance. This also means that as apparent from the proof of Theorem~\ref{thm:TrainingMIQP}, unlike usual in the big-$M$ formulation, the $M$ does not influence the feasibility of the problem. However, the value of $M$ influences the approximation error if it is not sufficiently large. In fact, a small $M$ leads to smaller values in the entries of the weights $\Vb$ and $\Wb$. Thus, $M$ implicitly regularizes the weights. Assuming the PWA case, we have the following result.
\begin{theorem}\label{thm:boundedVW}
    Let $\gb(\xb)$ and $\hb(\xb)$ be as in \eqref{eq:ghForPWA}. Then the norm of the matrices $\Vb_{k,1:n}$ and $\Wb_{l,1:n}$ of the OP~\eqref{eq:OP_training_MIQP} are bounded by the constants $v_k$ and $w_l$, respectively, i.e., $\lVert \Vb_{k,1:n} \rVert_2 \leq v_k, \ \forall k \in \{1,\dots,p_1\}$ and $\lVert \Wb_{l,1:n} \rVert_2 \leq w_l, \ \forall l \in \{1,\dots,p_2\}$.
\end{theorem}
\begin{proof}
Assume that $\deltab^{(i)}_{k_i}=1$ holds for $\xb^{(i)}$, then we have according to \eqref{eq:ViGeqVki}
\vspace{-1mm}
\begin{equation}
    \label{eq:VLeqM}
    \Vb_{k_i,1:n} \xb^{(i)} + \Vb_{k_i,n+1} \leq \Vb_{k,1:n} \xb^{(i)} + \Vb_{k,n+1} + M, \quad \forall i \in \{1,\dots,N_D\}, \ \forall k \in \{1,\dots,p_1\}.
    \vspace{-1mm}
\end{equation}
We now construct vectors $\tilde{\xb}^{(1)},\dots,\tilde{\xb}^{(n)}$ such that they are an orthonormal basis containing the vector $\tilde{\xb}^{(\kappa)}:=\frac{\Vb_{k_i,1:n}^\top}{\lVert\Vb_{k_i,1:n}\rVert_2}$. By representing $\xb^{(i)}$ in the orthonormal basis we can rewrite \eqref{eq:VLeqM} as follows
\vspace{-1mm}
\begin{align*}
    \Vb_{k_i,1:n} \sum_{l=1}^n \tilde{\xb}^{(l)} \tilde{\xb}^{(l)\top} \xb^{(i)} + \Vb_{k_i,n+1} &\leq \Vb_{k,1:n} \xb^{(i)} + \Vb_{k,n+1}\!+\!M \\
    \Leftrightarrow \!\! \lVert \Vb_{k_i,1:n} \rVert_2  \lVert \tilde{\xb}^{(\kappa)} \rVert_2 \tilde{\xb}^{(\kappa)\top} \xb^{(i)} \!+\!\Vb_{k_i,1:n}\!\!\!\! \sum_{l=1,l\neq \kappa}^n\!\!\!\! \tilde{\xb}^{(l)} \tilde{\xb}^{(l)\top} \xb^{(i)}\!+\!\Vb_{k_i,n+1} &\leq \Vb_{k,1:n} \xb^{(i)} + \Vb_{k,n+1}\!+\!M
\end{align*}
\vspace{-1mm}
\begin{equation*}
    \label{eq:boundedV}
    \Leftrightarrow\!\!\lVert \Vb_{k_i,1:n} \rVert_2 \!\leq\! \frac{1}{\tilde{\xb}^{(\kappa)\top} \xb^{(i)}} \! \left(\!M\!+\!\Vb_{k,n+1}\!-\!\Vb_{k_i,n+1}\!+\!\big(\Vb_{k,1:n}\!-\!\Vb_{k_i,1:n}\!\!\sum\limits_{l=1,l\neq \kappa}^n\tilde{\xb}^{(l)} \tilde{\xb}^{(l)\top}\big) \xb^{(i)}\!\right)\!\!=:\!v_k
    \vspace{-1mm}
\end{equation*}
for all $i \in \{1,\dots,N_D\}$ and for all $k \in \{1,\dots,p_1\}$. The upper bound $w_l$ for $\lVert \Wb_{l,1:n} \rVert_2$ can be derived in a similar way, which completes the proof.   
\end{proof}
\vspace{-2mm}
Theorem~\ref{thm:boundedVW} shows that it is not necessary to explicitly include a regularization in the cost function as in \citet{Siahkamari2020} to ensure small weights since $M$ implicitly regularizes the weights. Figure~\ref{fig:PWA} shows the effect of this implicit regularization. For the example, the MIQP~\eqref{eq:OP_training_MIQP} is solved with $p_1=p_2=2$ and \eqref{eq:ghForPWA}, for different values of $M$. In Figure~\ref{fig:PWA}~(a), the constant is not sufficiently large and thus an exact fit of the data is impossible. For $M=50$ and $M=100$, an exact fit is possible. However, $M=50$ leads to a fit with smaller weights. 

Another advantage of the proposed method is that since for $\gb(\xb)$ and $\hb(\xb)$ as in \eqref{eq:ghForPWA}, \eqref{eq:Phi} is a maxout NN with affine pre- and post-activation, we can directly use \citet[Lem.~2]{Teichrib2023} to show that \eqref{eq:Phi} can be modelled as an MI feasibility problem and is thus well-suited for the use inside an OP.   
\begin{corollary}\label{Cor:MI}
    Let $\gb(\xb)$ and $\hb(\xb)$ be as in \eqref{eq:ghForPWA}. Then, the function \eqref{eq:Phi} can be modelled as an MI feasibility problem. 
\end{corollary}
\begin{proof}
    According to \citet[Lem.~2]{Teichrib2023}, an MI feasibility problem with the constraints \citet[Eq.~(18)]{Teichrib2023} can be used to replace the nonlinear relation \eqref{eq:Phi}. 
\end{proof}
\vspace{-2mm}
This corollary allows the replacement of the nonlinear function \eqref{eq:Phi} inside an OP by MI linear constraints. Thus, a nonlinear OCP, where \eqref{eq:Phi} is used either as a prediction model or to approximate the cost function, can be reformulated as MIP. This enables the use of standard MI solvers to find a global optimum of the OP.

\vspace{-2mm} 
\section{Complexity reduction}

In this section, we describe two approaches to reduce the complexity of the MIQP~\eqref{eq:OP_training_MIQP}. The approaches mainly focus on reducing the number of binary variables. Before describing our approaches, we first show that the MIQP~\eqref{eq:OP_training_MIQP} can be solved by a QP without using any binary variables if we are not interested in a small number of segments $p_1$ and $p_2$ in \eqref{eq:Phi}. Consider the special case where we have exactly as many segments as data points, i.e. $p_1=p_2=N_D$. In this case, the QP
\vspace{-2mm}
\begin{subequations}\label{eq:OP_training_QP}
    \begin{align}
    \min_{\Vb,\Wb,\alpha^{(i)},\beta^{(i)}} \sum_{i=1}^{N_D}& \big(y^{(i)}-(\alpha^{(i)}-\beta^{(i)})\big)^2 \\
    \label{eq:VWLeq1QP}
    \text{s.t.} \quad \Vb \gb(\xb^{(i)}) &\leq \oneb \alpha^{(i)}, \hspace{4.9mm} \Wb \hb(\xb^{(i)}) \leq \oneb \beta^{(i)}, \quad \forall i \in \{1,\dots N_D\}, \\
    \label{eq:VWKappa}
    \Vb_{\kappa} \gb(\xb^{(\kappa)}) &\geq \alpha^{(\kappa)}, \quad \Wb_{\kappa} \hb(\xb^{(\kappa)}) \geq \beta^{(\kappa)}, \hspace{4.5mm} \forall \kappa \in \{1,\dots,N_D\}
    \vspace{-2mm}
    \end{align}
\end{subequations}
can be used to compute optimal $\Vb$, $\Wb$, $\alpha^{(i)}$, and $\beta^{(i)}$ of the MIQP~\eqref{eq:OP_training_MIQP}. Since there is one segment for each data point, no binary variables are needed to decide which segment is active, thus \eqref{eq:VgLeqAlpha}-\eqref{eq:sumGamma} can be replaced by \eqref{eq:VWLeq1QP} and \eqref{eq:VWKappa}. However, this leads to the fact that the number of data points determines the number of segments. As a result, we have $2 N_D(N_D+1)$ inequality constraints and $2N_D + N_D (r_1 + r_2)$ optimization variables in \eqref{eq:OP_training_QP}. The QP~\eqref{eq:OP_training_QP} can be reformulated to obtain the QP from \citet{Siahkamari2020}. Thus, our method generalizes known results for piecewise linear regression and includes the OP from \citet{Siahkamari2020} as a special case. The elimination of the binary variables in the OP~\eqref{eq:OP_training_QP} is only possible for $p_1=p_2=N_D$, which leads to a complex function $\Phi(\xb)$ with a high number of segments, for large data sets. Therefore, in the remainder of this section, we will present two approaches for reducing the number of binary variables in the MIQP~\eqref{eq:OP_training_MIQP} without increasing the complexity of the fitted function $\Phi(\xb)$, in terms of the number of segments. 

\subsection{Symmetry breaking}\label{sec:SymBreak}
Due to symmetries in the MIQP~\eqref{eq:OP_training_MIQP}, there are many equivalent solutions. These symmetries result from the fact that the values of the max-functions in \eqref{eq:Phi} are independent of the order of the $p_1$ and $p_2$, respectively, segments. In the following proposition, we will present a method to eliminate some of these redundant solutions.
\begin{proposition}\label{prop:SymBreak}
    There exists an optimal solution of the MIQP~\eqref{eq:OP_training_MIQP} with
    \begin{equation}\label{eq:SymBreak}
    \deltab_{2:p_1}^{(1)}=\zerob, \ \deltab_{3:p_1}^{(2)}=\zerob, \ \dots, \ \deltab_{p_1}^{(p_1-1)}=0 \quad \text{and} \quad \gammab_{2:p_2}^{(1)}=\zerob, \ \gammab_{3:p_2}^{(2)}=\zerob, \ \dots, \ \gammab_{p_2}^{(p_2-1)}=0.
    \vspace{-1mm}
    \end{equation}
\end{proposition}
\begin{proof}
    Every permutation of the $k$-indices of an optimal solution of the MIQP~\eqref{eq:OP_training_MIQP} is also an optimal solution since the maximum is invariant under a permutation of the segments. Assume $\deltab^{(i)}$ and $\gammab^{(i)}$ are the binary variables of an optimal solution of the MIQP~\eqref{eq:OP_training_MIQP}. Due to the discussed invariance of the maximum, we can assume that the optimal solution is such that the first segment is active for the first data point, i.e., $\deltab_1^{(1)}=1$ and $\gammab_1^{(1)}=1$. In combination with \eqref{eq:sumDelta} and \eqref{eq:sumGamma} this leads to $\deltab_{2:p_1}^{(1)}=\zerob$ and $\gammab_{2:p_2}^{(1)}=\zerob$. For the second data point, there are two possibilities. Either the same segment is active as for the first data point, or a different segment is active, which we can assume to be the second segment, i.e., $\deltab_1^{(2)}+\deltab_2^{(2)}=1$ and $\gammab_1^{(2)}+\gammab_2^{(2)}=1$ or equivalent $\deltab_{3:p_1}^{(2)}=\zerob$ and $\gammab_{3:p_2}^{(2)}=\zerob$. We can proceed in a similar way with the remaining data points until the $(p_1-1)$-th and $(p_2-1)$-th, respectively, data point, which leads to the remaining constraints in \eqref{eq:SymBreak} and completes the proof. 
\end{proof}
\vspace{-1mm}
According to Proposition~\ref{prop:SymBreak} we can reduce the number of binary variables by including the constraints \eqref{eq:SymBreak} in the MIQP~\eqref{eq:OP_training_MIQP}. Symmetry breaking is also achieved in \citet[Prop.~3.2]{Toriello2012} by enforcing an order to the parameter of the fitted function through additional constraints. However, unlike our method, this does not reduce the number of binary variables.

\vspace{-1mm}
\subsection{Preclustering}\label{sec:PreClustering}

Since the number of data points determines the number of binary variables in \eqref{eq:OP_training_MIQP}, reducing the number of data points is a very effective way to reduce the computation time needed to solve the MIQP~\eqref{eq:OP_training_MIQP}. Unfortunately, a certain number of data points is often required to obtain an approximation that generalizes well to unseen data. This means that it is impossible to simply reduce the size of the data set. However, we can precompute a clustering of the data set and modify the OP so that instead of assigning individual data points to segments of the max-function, whole clusters of multiple data points are assigned. Assume the data set is partitioned into $K$ disjoint clusters of the form $\Cc_j:=\{ (\xb^{(i)}, y^{(i)}) \in \Dc \ | \ i \in \Ic_j \}$ with $\Ic_j\subseteq\{1,\dots,N_D\}$, $\Cc_i\cap\Cc_j=\emptyset$ for all $i \neq j$, and $\cup_{j=1}^K \Cc_j = \Dc$. Then, we can formulate a modified version of MIQP~\eqref{eq:OP_training_MIQP} with    
\vspace{-1mm}
\begin{subequations}\label{eq:PreClustering}
\begin{align}
    \label{eq:deltaCluster}
    \Vb \gb(\xb^{(i)}) &\leq \oneb \alpha^{(i)}, \hspace{5.8mm} \Vb \gb(\xb^{(i)}) + M (\oneb-\deltab^{(l)}) \geq \oneb \alpha^{(i)}, \quad \forall i \in \{1,\dots N_D\}, \\
    \label{eq:gammaCluster}
    \Wb \hb(\xb^{(i)}) &\leq \oneb \beta^{(i)}, \quad \Wb \hb(\xb^{(i)}) + M (\oneb-\gammab^{(l)}) \geq \oneb \beta^{(i)}, \hspace{4.1mm} \forall i \in \{1,\dots N_D\},
    \vspace{-1mm}
    \end{align}
\end{subequations}
instead of \eqref{eq:VgLeqAlpha} and \eqref{eq:WhLeqBeta}, where $l$ is such that we have $(\xb^{(i)}, y^{(i)}) \in \Cc_l$. The binary variable vectors $\deltab^{(l)} \in \{0,1\}^{p_1}$ and $\gammab^{(l)} \in \{0,1\}^{p_2}$ select which segment of the max-functions is active for the whole cluster $\Cc_l$ and not just for a single data point as in the initial MIQP~\eqref{eq:OP_training_MIQP}. Any algorithm that computes a clustering with the desired properties can be used here, e.g. k-means$++$ from \citet{Arthur2007}. Note that although we perform a clustering before solving the optimization problem, the modified MIQP~\eqref{eq:OP_training_MIQP} with \eqref{eq:PreClustering} instead of \eqref{eq:VgLeqAlpha} and \eqref{eq:WhLeqBeta} does not computes a linear regression for each cluster individually as in \citet{Magnani2009}, since there is not a distinct segment for each cluster. As long as $p_i<K$ with $i \in \{1,2\}$, multiple clusters are assigned to the same segment.
With the preclustering, we can reduce the number of binary variables from $N_D(p_1+p_2)$ to $K(p_1+p_2)$. For $K=N_D$, i.e., each cluster contains exactly one data point, the constraints \eqref{eq:VgLeqAlpha}, \eqref{eq:WhLeqBeta} and \eqref{eq:deltaCluster}, \eqref{eq:gammaCluster} are equivalent. Furthermore, analogous to \eqref{eq:OP_training_QP} for $p_1=p_2=K$ the MIQP with preclustering can be reformulated as a QP by replacing \eqref{eq:VWKappa} in the QP~\eqref{eq:OP_training_QP} with
\vspace{-1mm}
\begin{equation}
    \label{eq:activeSegQPCluster}
    \Vb_{\kappa} \gb(\xb^{(i)}) \geq \alpha^{(i)}, \quad \Wb_{\kappa} \hb(\xb^{(i)}) \geq \beta^{(i)}, \quad \forall i \in \{1,\dots N_D\},
    \vspace{-1mm}
\end{equation}
where $\kappa$ in \eqref{eq:activeSegQPCluster} is such that we have $(\xb^{(i)}, y^{(i)}) \in \Cc_{\kappa}$. The resulting QP has $2 N_D (K + 1)$ inequality constraints and $2 N_D + K (r_1 + r_2)$ optimization variables. 

\vspace{-2mm}
\section{Numerical examples}
In this section, we use the MIQP~\eqref{eq:OP_training_MIQP} with different $\gb(\xb)$ and $\hb(\xb)$ to approximate a nonlinear MPC and compare the results with an approximation by an NN.

\subsection{Nonlinear MPC}
\begin{figure}[h!]
\centering
    \subfigure[\hspace{0mm} \label{fig:1}]{\includegraphics[trim={5cm 11cm 5cm 11cm},clip,scale=0.45]{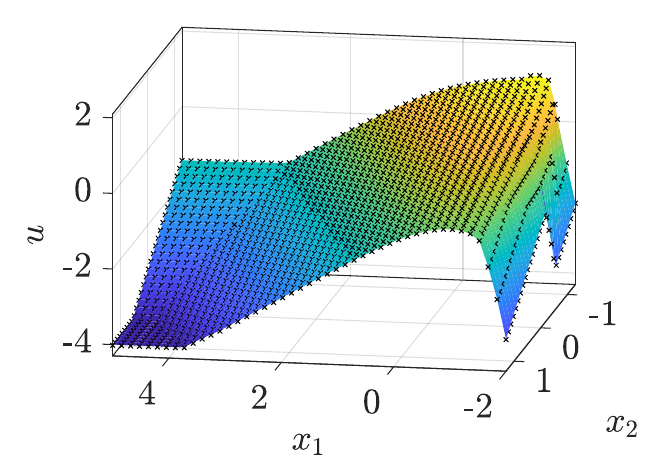}}
    \subfigure[\hspace{0mm}]{\includegraphics[trim={5cm 11cm 5cm 11cm},clip,scale=0.45]{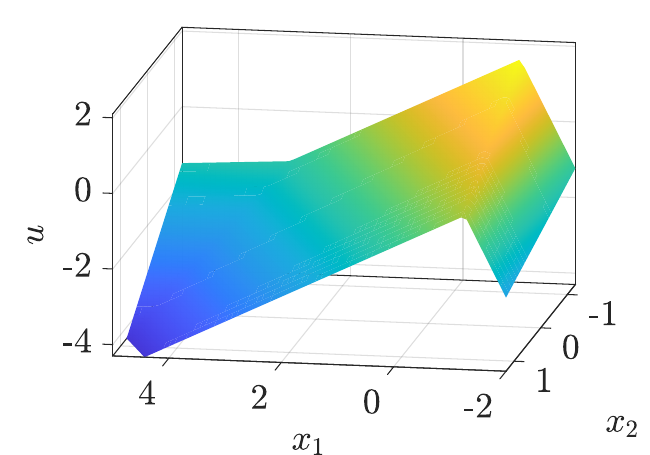}}
    \subfigure[\hspace{0mm}]{\includegraphics[trim={5cm 11cm 5cm 11cm},clip,scale=0.45]{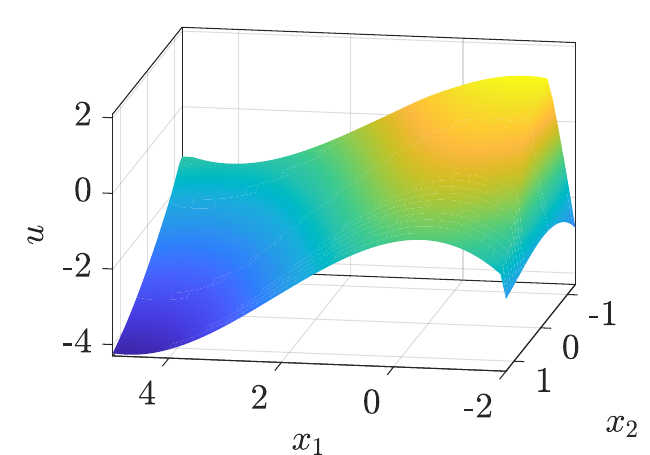}}
    \vspace{-2mm}
\caption{(a) shows the optimal control law with the sampled data points. (b) and (c) are approximations of the optimal control law by the functions $1$ and $4$, respectively, from Table~\ref{tab:exampleApprox}.}\label{fig:OCL}
\vspace{-4mm}
\end{figure}

Consider the nonlinear OCP 
\vspace{-2mm}
\begin{align}
    \label{eq:OCP}
    V(\xb) := \!\!\!\!\min_{\substack{\xb(0),\dots,\xb(N)\\ u(0),\dots,u(N-1)}} 
    \!\!\!\!\!\!\!\!\!\!\!\!\!\!\!\!\!\!\!\! & 
    \,\,\,\,\,\,\,\,\,\,\,\,\,\,\,\,\,\,\,\,\,
    \xb(N)^\top\begin{pmatrix} 1.3 & 1.9 \\ 1.9 & 3.0 \end{pmatrix}\xb(N) + \! \sum_{k=0}^{N-1} \lVert\xb(k)\rVert_2^2 + \lVert u(k)\rVert_2^2   \\
    \nonumber
    \text{s.t.} \ \xb_1(k+1)&=\xb_1(k) + T_s \xb_2(k), \hspace{58mm} \forall k \in \{0,\dots,N-1\}, \\
    \nonumber
    \xb_2(k+1)&=\xb_2(k) - T_s \frac{\rho_0}{m}\mathrm{e}^{-\xb_1(k)}\xb_1(k)-T_s\frac{h_d}{m}\xb_2(k)+T_s\frac{u(k)}{m}, \quad\!\!\forall k \in \{0,\dots,N-1\}, \\
    \nonumber
    -1.3&\leq\xb_2(k)\leq1.3, \quad -4\leq u(k)\leq 4, \hspace{35.7mm} \forall k \in \{0,\dots,N-1\}, \\
    \nonumber
    \xb(0)&=\xb, \quad \xb(N)^\top\begin{pmatrix} 1.3 & 1.9 \\ 1.9 & 3.0 \end{pmatrix}  \xb(N) \leq 0.001,
    \vspace{-2mm}
\end{align}
from \citet{Magni2003} with parameters $N=15$, $T_s=0.4$, $m=1$, $\rho_0=0.33$ and $h_d=1.1$. The OCP~\eqref{eq:OCP} is solved for $1350$ different values of $\xb$ on a regular grid of size $45 \times 30$ with $\xb_1$ from $-2$ to $5$ and $\xb_2$ from $-1.3$ to $1.3$ to generate a data set of the form $\big(\xb^{(i)},y^{(i)}=u^\ast(\xb^{(i)})\big)$, where $u^\ast(\xb^{(i)})$ is the optimal $u(0)$ for the state $\xb=\xb^{(i)}$. The nonlinear OCP~\eqref{eq:OCP} is solved using the software framework CasADi from \citet{Andersson2019}. Based on the data set, two different maxout NN of the form \eqref{eq:Phi} with $\gb(\xb)$ and $\hb(\xb)$ as in \eqref{eq:ghForPWA} are trained using the Adam algorithm \citep{kingma2017adam}. The values of $p_1$ and $p_2$ of the two NN are chosen as $p_1=p_2=2$ and $p_1=p_2=50$ (row $1$ and $2$ of Table~\ref{tab:exampleApprox}). Such NN can be interpreted as maxout NN with affine pre- and post-activation and $2$ neurons of rank $2$ and $50$, respectively, in one hidden layer. We also used the MIQP~\eqref{eq:OP_training_MIQP} with $p_1=p_2=2$, different choices of $\gb(\xb)$ and $\hb(\xb)$ and the symmetry breaking constraints \eqref{eq:SymBreak} to compute global optimal parameters $\Vb$ and $\Wb$ of the function $\Phi(\xb)$. We start with $\gb(\xb)$ and $\hb(\xb)$ as in \eqref{eq:ghForPWA}, which results in a function that is equal to the first maxout NN, and add monomials up to the order of three to the functions (row $3$ and $4$ of Table~\ref{tab:exampleApprox}). The results are summarized in Table~\ref{tab:exampleApprox}. The last two columns show the quality of the approximation in terms of the mean squared error (MSE) and the maximum absolute value of the error on a test data set. The second row indicates the method used to compute the parameters.
\begin{table*}[h]
    \vspace{-4mm}
    \caption{Different approximations of the nonlinear optimal control law $u^\ast(\xb)$.}
    \vspace{2mm}
    \centering
    \begin{tabular}{ccclcc}
        \toprule
        \!\!\!\!No.\!\!\!\! & Method & $p_1,p_2$ & $\gb(\xb),\hb(\xb)$ & MSE & $e_{\text{max}}$ \\  
        \midrule
        \!\!\!\!$1.$\!\!\!\! & NN & \!\!\!\!$2$\!\!\!\! & $\begin{pmatrix} 1 & \xb^\top \end{pmatrix}^\top$ & \!\!$0.0532$\!\! & $1.6054$ \\[.5ex]
        \!\!\!\!$2.$\!\!\!\! & NN & \!\!\!\!$50$\!\!\!\! & $\begin{pmatrix} 1 & \xb^\top \end{pmatrix}^\top$ & \!\!$0.0252$\!\! & $1.4971$ \\[.5ex]
        \!\!\!\!$3.$\!\!\!\! & MIQP~\eqref{eq:OP_training_MIQP} & \!\!\!\!$2$\!\!\!\! & $\begin{pmatrix} 1 & \xb^\top \end{pmatrix}^\top$ & \!\!$0.0503$\!\! & $1.0797$ \\[.5ex]
        \!\!\!\!$4.$\!\!\!\! & MIQP~\eqref{eq:OP_training_MIQP} & \!\!\!\!$2$\!\!\!\! & $\begin{pmatrix} 1 & \xb^\top & \xb_1^2 & \xb_1 \xb_2 & \xb_2^2 & \xb_1^3 & \xb^2_1 \xb_2 & \xb_1 \xb^2_2 & \xb^3_2 \end{pmatrix}^\top$ & \!\!$0.0240$\!\! & $1.0687$ \\[.5ex]
        \bottomrule
        \label{tab:exampleApprox}
    \end{tabular}
    \vspace{-4mm}
\end{table*}
As apparent from Table~\ref{tab:exampleApprox}, the approximation provided by the MIQP~\eqref{eq:OP_training_MIQP} with cubic terms in $\gb(\xb)$ and $\hb(\xb)$ gives the lowest error, even though it has only two segments. By increasing the number of segments $p_1$ and $p_2$, the MSE of the first NN approaches the MSE of the $4$th function. However, the maximum error barely changes, which may be a problem since it is shown in \citet{Fabiani2022} that the maximum error is crucial for the stability and performance of an approximated controller. Another advantage of the proposed method is that adding additional terms to the $\gb(\xb)$ and $\hb(\xb)$ functions can improve the approximation quality without increasing the number of segments. This is beneficial for a fast evaluation of the approximation $\Phi(\xb)$. Figure~\ref{fig:OCL} shows how the function $4$ from Table~\ref{tab:exampleApprox} leads to a better result by approximating the curvature of the optimal control law.

\vspace{-4mm}
\section{Conclusion}
We proposed a new MIP-based method for piecewise regression, which allows to fit an arbitrary piecewise defined function to data points by using the difference of two max-functions \eqref{eq:Phi}. The presented method is very flexible as it is not restricted to a special class of functions and thus allows to trade off the complexity of the fitted function against the approximation error. In addition, we showed in a numerical example that the method can be used to find a simple approximation, in terms of the number of regions, of the optimal control law of nonlinear MPC.   



\end{document}